\documentclass[a4paper,12pt,reqno]{amsart}
\usepackage{amsmath, amssymb,amsthm}
\usepackage{enumerate}
\usepackage{cite}
\usepackage{graphicx}
\usepackage{subcaption}
\numberwithin{equation}{section}
\newtheorem{theorem}{Theorem}[section]
\newtheorem{lemma}[theorem]{Lemma}

\newtheorem{remark}[theorem]{Remark}
\newtheorem{definition}[theorem]{Definition}
\theoremstyle{definition}

\begin{document}
\title[Self-inversive polynomials and quasi-orthogonality on the unit circle]
{Self-inversive polynomials and quasi-orthogonality on the unit circle}
\author[K. K. Behera]{Kiran Kumar Behera}
\address{
Department of Mathematics,
Indian Institute of Science Bangalore, India}
\email{kiranbehera@iisc.ac.in}
\thanks{
This research is supported by the Dr. D. S. Kothari postdoctoral
fellowship scheme of University Grants Commission (UGC), India.}
%
%
\begin{abstract}
In this paper we study quasi-orthogonality on the unit circle based on the structural and orthogonal properties of a class of self-invariant polynomials.
We discuss a special case in which these polynomials are represented in terms of the reversed Szeg\H{o} polynomials 
of consecutive degrees and illustrate the results using contiguous relations of hypergeometric functions.
This work is motivated partly by the fact that
recently cases have been made to establish para-orthogonal polynomials as the unit circle analogues of quasi-orthogonal polynomials on the real line so far as spectral properties are concerned.
We show that structure wise too there is great analogy when self-inversive polynomials are used to study quasi-orthogonality on the unit circle.
\end{abstract}
\maketitle
\section{Introduction}
The concept of quasi-orthogonality of polynomials was introduced by Riesz and later studied, among others, by Fej\`{e}r, Shohat and Chihara in relation to moment problems and associated quadrature formulae.
If $\{P_n(x)\}_{n=0}^{\infty}$ is a sequence of polynomials orthogonal with respect to a positive weight function $w(x)$
on the real line $\mathbb{R}$,
then a necessary and sufficient condition
\cite{Brezinski-Driver-QOP-2004}
for a polynomial $Q_n(x)$ to be quasi-orthogonal of order $r$ with respect to $w(x)$ is 

\begin{align}
	\label{eqn: quasi-orthogonal linear combination}
	Q_n(x)=c_{n,0}P_n(x)+c_{n,1}P_{n-1}(x)+\cdots+
	c_{n,r}P_{n-r}(x),
	\quad n\geq r,	
\end{align}
where $c_{n,i}\in\mathbb{R}$ and $c_{n,0}c_{n,r}\neq0$.
The relation \eqref{eqn: quasi-orthogonal linear combination}
justifies the name quasi-orthogonality: $Q_n(x)$ is orthogonal to every polynomial of degree not exceeding $n-r$
with respect to $w(x)$. This can be equivalently stated as
\begin{align}
	\label{eqn: quasi-orthogonality integral}
	\int_{a}^{b}x^kQ_n(x)w(x)dx=
	\begin{cases} 
		0 &\mbox{if } k=0,\cdots,n-r-1, \\
		h_n\neq0 & \mbox{if } k=n-r, 
	\end{cases}
\end{align}
where we have the standard orthogonality on 
$\mathbb{R}$ for $r=0$. 
We refer to
\cite{Brezinski-Driver-QOP-2004,
Chihara-QOP-PAMS-1957,
Dickinson-QOP-PAMS-1961,
Draux-quasi-orthogonality-JAT-1990,
Ismail-Wang-QOP-ODE-2019,
Joulak-QOP-ANM-2005,
Shohat-QOP-1937,
Xu-QOP-JMAA-1994}
and references therein for this classical theory and its applications.

There have been attempts to generalize the concept of
quasi-orthogonality from 
$\mathbb{R}$ to the unit circle $\partial{\mathbb{D}}$
via the relations 
\eqref{eqn: quasi-orthogonal linear combination}
and \eqref{eqn: quasi-orthogonality integral}.
However, it was proved that if
$\{P_n(z)\}$ is a sequence of polynomials orthogonal on the unit circle, also called Szeg\H{o} polynomials
\cite{Simon-book-Part-1,Szego-book},
then the necessary and sufficient condition 
\cite[Theorem 1]{Marcellan-QOP-Berstein-1996}
(see also
\cite{Alfaro-Moral-1994,
	Branquinho-Marcellan-QOP-1996})
for the polynomial $Q_n(z)$ given by 
\eqref{eqn: quasi-orthogonal linear combination}
to be orthogonal with respect to a 
non-trivial positive measure on
the unit circle is that $Q_n(z)$ 
belongs to the class of Bernstein-Szeg\H{o} polynomials.
Hence, a new concept of quasi-orthogonality is defined in
\cite{Alfaro-Moral-1994} which is dependent on the structure 
of the Szeg\H{o} polynomials and related semi-classical forms are discussed.

The Szeg\H{o} polynomials find many applications in areas like
approximation on the complex plane $\mathbb{C}$, prediction theory and signal processing
\cite{
DG-OPUC-signal-1989,
Dewilde-Dym-lossless-prediction-1981,
Grenander-Szego-book-Toeplitz-form,
Jones-Njastad-Szego-digital-1991,
Jones-Njastad-Waadeland-Szego-frequency-1994
}. 
However, they suffer from the major drawback that their zeros lie outside the support of the measure and hence cannot be used in interpolation processes on the unit circle. 
This was overcome by introducing para-orthogonal polynomials which are self-invariant polynomials 
with symmetric orthogonality conditions.
The polynomial $\mathcal{P}_n(z)$ is a self-invariant polynomial if and only if
it satisfies $\mathcal{P}_n^{\ast}(z)=\tau_n\mathcal{P}_n(z)$,
where $\tau_n\in\partial{\mathbb{D}}$ and $\mathcal{P}_n^{\ast}(z)=
z^n\overline{\mathcal{P}_n(1/\bar{z})}$.
This invariance property leads to the symmetry that 
$\mathcal{P}_n(z)$ is orthogonal to the monomial $z^k$ if and only if $\mathcal{P}_n(z)$ is orthogonal to the monomial $z^{n-k}$ for $k=0,\cdots,n$.
Hence, if one defines the spaces 
\begin{align*}
\Lambda_{n,2l+1}=\mbox{span}\{z^k:k=l+1,\cdots,n-l-1\},
\quad
0\leq l\leq \left\lfloor\frac{n}{2}\right\rfloor-1,
\end{align*}
then the invariant polynomials lying in 
$\Lambda_{n,1}^{\perp}$ $(l=0)$ are precisely the 
para-orthogonal polynomials.
These are introduced in
\cite{Jones-Njastad-Thron-moment-BLMS-1989} to
study quadrature formula on the unit circle
and in the solution of the trigonometric moment problem.
They are also used to adapt many spectral properties
developed on the real line to the unit circle
\cite{CMV-measures-POP-2002,
	Breuer-POP-singular-JAT-2020,
	Golinskii-POP-2002} 

In view of this, cases have been made, see for instance
\cite{Bultheel-quasi-POP_JCAM-2022}, to establish
para-orthogonal polynomials as the counterpart of
quasi-orthogonal polynomials on $\mathbb{R}$ so far as
spectral properties are concerned.
In particular, the classical concept of 
para-orthogonality is generalized 
\cite{Bultheel-quasi-POP_JCAM-2022} 
to the concept of 
quasi para-orthogonal polynomials
of order $2l+1$, which is precisely 
the set of invariant polynomials of 
degree $n$ lying in the space 
$\Lambda_{n,2l+1}^{\perp}$.
The prefix quasi refers to the fact that
for fixed degree $n$, these polynomials satisfy fewer 
(symmetric) orthogonality conditions.

Motivated by the
assertion that invariance of polynomials 
should have a role in the unit circle analogue of quasi-orthogonality 
on $\mathbb{R}$,
we study quasi-orthogonality 
on $\partial{\mathbb{D}}$ based on the structure of a class of
self-invariant polynomials.
Such polynomials satisfy a three term recurrence relation of the form
\begin{align}
\label{eqn: ttrr for DG-invariant polynomials}
\tilde{R}_{n+1}(z)=(\bar{\beta}_nz+\beta_n)\tilde{R}_{n}(z)
-z\tilde{R}_{n-1}(z),
\quad
n\geq0,
\end{align}
and are introduced by Delsarte and Genin 
\cite{DG-OPUC-signal-1989,
DG-split-Levinson-1986,
DG-tridiagonal-algebraic-I-1991,
DG-tridiagonal-Szego-1988,
DG-trdiagonal-algebraic-II-1991}
in order to solve certain 
problems in digital signal processing.
Hence, we will refer to this class as the 
$\mathcal{DG}$ 
class of invariant polynomials.
This class is also used to 
characterize a family of non-trivial measures on
the unit circle along with
the corresponding families of para-orthogonal polynomials 
and the Szeg\H{o} polynomials
\cite{Ranga-Favard-type=JAT-2014,
Ranga-OPUC-chain-seq-JAT-2013}.

The key to our study is the following representation 
\begin{align}
\label{eqn: introduction-linear-combination}
(z-1)\mathcal{P}_n(z)=\varphi_{n+1}^{\ast}(z)+
c_{n+1}\varphi_{n}^{\ast}(z),
\quad
n\geq1,
\end{align}
that we obtain for any polynomial $\mathcal{P}_n(z)$ in  
$\mathcal{DG}$ class, where $\varphi_n^{\ast}(z)$
is the monic reversed Szeg\H{o} polynomial
of degree $n$.
The crucial result we prove in this context is 
the following.

\emph
{	
Given any complex parameter $\zeta$ lying on the unit circle,
there exists a polynomial sequence $\{\mathcal{B}_n(z)\}_{n=0}^{\infty}$
such that the polynomials defined by
\begin{align}
\label{eqn: explicit representation in introduction}
\mathcal{P}_0^{(\zeta)}(z):=1,
\quad
\mathcal{P}_n^{(\zeta)}(z)=\frac{1}{z-\zeta}
\left(
\mathcal{B}_{n+1}(z)+
\omega_{n+1}(\zeta)\mathcal{B}_n(z)
\right),
\quad
n\geq1,
\end{align}
for some $\omega_{n}(\zeta)\in\mathbb{C}$ 
belong to the $\mathcal{DG}$ class of invariant polynomials.
}

The polynomials $\mathcal{P}_n(z)$
given by \eqref{eqn: introduction-linear-combination} 
are obtained from 
$\mathcal{P}_n^{(\zeta)}(z)$ for $\zeta=1$.
Our use of $\mathcal{DG}$ class serves two purposes. First,
the representation \eqref{eqn: introduction-linear-combination}
is similar to \eqref{eqn: quasi-orthogonal linear combination}
for $r=1$ and hence $\mathcal{P}_n(z)$ is at least structurally similar to a quasi-orthogonal polynomial of order $1$ on the real line. 
Second, the use of symmetric orthogonality conditions 
relaxes the quasi-orthogonality conditions on the unit circle.
Further, we show that $\mathcal{P}_n^{(\zeta)}(z)$, $n\geq1$, satisfy a three term recurrence relation,
a result that is known for a quasi-orthogonal polynomial of order $1$ on the real line
\cite{Chihara-QOP-PAMS-1957,
	Draux-quasi-orthogonality-JAT-1990}.
In addition, we prove that if 
$\mathcal{P}_n(z)$ is quasi-orthogonal with respect to
$\mathfrak{v}$ and satisfy symmetric orthogonality conditions with respect to $\mathfrak{u}$, then
\begin{align}
	\label{eqn: characterize v in terms of u introduction}
	\mathfrak{v}=[z\mathcal{U}(z)+z^{-1}
	\bar{\mathcal{U}}(z^{-1})]\mathfrak{u},
	\quad
	\mathcal{U}(z)=
	\frac{1}{2}(c_{s-1}z^s+\cdots+c_1z^2+c_0z+\bar{c}_1),
\end{align}
where $c_{s-1}\neq0$ so that $\mathcal{U}(z)$ is a polynomial of degree $s$.

The manuscript is organized as follows. 
In Section~\ref{sec: Structural relations}, we obtain preliminary structural relations that lead to the proof
of 
\eqref{eqn: explicit representation in introduction}.
In particular, we give an algorithm to generate the sequence
$\omega_n(\zeta)$, $n\geq2$, with a special choice of expression for $\omega_1(\zeta)$. 
In Section~\ref{sec: special case},
we identify that the class of reversed Szeg\H{o} polynomials can be used for $\mathcal{B}_n(z)$ in the representation 
\eqref{eqn: explicit representation in introduction}. 
We prove a Szeg\H{o} type relation that any polynomial in the $\mathcal{DG}$ class satisfies.
We state our definition of quasi-orthogonality on $\partial{\mathbb{D}}$ in 
Section~\ref{sec: quasi-orthogonality}
and use this to obtain the characterization 
\eqref{eqn: characterize v in terms of u introduction}. 
We conclude with Section~\ref{sec: Illustration}
in which we illustrate the theory presented in the paper
using contiguous hypergeometric relations.
\section{Structural Relations}
\label{sec: Structural relations}
Let $\{\mathcal{B}_n(z)\}_{n=0}^{\infty}$ be a sequence of polynomials given by the three term recurrence relation
\begin{align}
\label{eqn: ttrr for Bn}
\mathcal{B}_{n+1}(z)=(z+\sigma_{n+1})\mathcal{B}_n(z)-
\lambda_{n+1}z\mathcal{B}_{n-1}(z),
\quad
n\geq1,
\end{align}
with $\mathcal{B}_{-1}(z)=0$ and $\mathcal{B}_0(z)=1$. 
Even though $\lambda_1$ does not affect 
\eqref{eqn: ttrr for Bn}, we fix $\lambda_1\in\mathbb{C}\setminus\{0\}$ and consider the complex parameters $\sigma_n$ and $\lambda_n$, $n\geq1$.
A characterization of such polynomials is that if $\sigma_n\neq0$ and $\lambda_{n+1}\neq0$, $n\geq1$, then there exists a quasi-definite moment functional $\mathcal{M}$ such that
\begin{align*}
\mathcal{M}[z^{-k}\cdot\mathcal{B}_n(z)]=\delta_{n,k}
\frac{\lambda_2\cdots\lambda_{n+1}}{\sigma_2\cdots\sigma_{n+1}},
\quad
k=0,1,\cdots,n,
\quad
n\geq1.
\end{align*}
Further, with specific choices for $\sigma_n$ and $\lambda_{n+1}$, the recurrence relation 
\eqref{eqn: ttrr for Bn} can be transformed into the form
\eqref{eqn: ttrr for DG-invariant polynomials}
satisfied by $\mathcal{B}_n(z)$ after a scaling.

Let $\zeta\in\mathbb{C}\setminus\{0\}$ be fixed
and $z\in\mathbb{C}\setminus\{\zeta\}$. 
We consider the polynomial sequence $\{\mathcal{Q}_n^{(\zeta)}(z)\}_{n=0}^{\infty}$ defined by
\begin{align}
\label{Qn as a linear combination of Bn}
\mathcal{Q}_0^{(\zeta)}(z):=1,
\quad
\mathcal{Q}_n^{(\zeta)}(z):=\mathcal{B}_n(z)+
\omega_n(\zeta)\mathcal{B}_{n-1}(z),
\quad
n\geq1,
\end{align}
where $\{\omega_n(\zeta)\}_{n=1}^{\infty}$ is a sequence of complex numbers depending on $\zeta$. 

The next result is partly motivated by
\cite{Draux-quasi-orthogonality-JAT-1990}, 
precisely the fact that quasi-orthogonal polynomials on 
$\mathbb{R}$ satisfy a three term recurrence relation with
polynomial coefficients.
\begin{lemma}
\label{lemma for mixed recurrence relations}
Suppose $\omega_{n}(\zeta)$, $n\geq2$,
is defined recursively as
\begin{align}
\label{recursive definition of omega-n}
\omega_{n}(\zeta)=
-\sigma_n-\zeta
\left(1+\frac{\lambda_n}{\omega_{n-1}(\zeta)}\right)
\quad n\geq2,
\end{align}
where $\omega_1(\zeta)$
is arbitrary.
Then
$\{\mathcal{Q}_{n+1}^{(\zeta)}(z)\}_{n=2}^{\infty}$
satisfy a three term recurrence relation of the form
\eqref{eqn: ttrr for Bn}
with polynomial coefficients.
\end{lemma}
\begin{proof}
For $n\geq1$, consider the system
\begin{align}
\label{system of equations for mixed ttrr}
\left(
\begin{array}{ccccc}
0 & 0 & 0 & \omega_{n+1}(\zeta) & 1 \\
0 & 0 & \omega_n(\zeta) & 1 & 0 \\
-1 & \omega_{n-1}(\zeta) & 1 & 0 & 0 \\
0 & 0 & -\tau_nz & z-\beta_n & -1 \\
0 & -\tau_{n-1}z & z-\beta_{n-1} & -1 & 0 \\
\end{array}
\right)
\left(
\begin{array}{c}
\mathcal{Q}_{n-1}^{(\zeta)}(z) \\
\mathcal{B}_{n-2}(z) \\
\mathcal{B}_{n-1}(z) \\
\mathcal{B}_{n}(z) \\
\mathcal{B}_{n+1}(z) \\
\end{array}
\right)
=
\left(
\begin{array}{c}
\mathcal{Q}_{n+1}^{(\zeta)}(z) \\
\mathcal{Q}_{n}^{(\zeta)}(z) \\
0 \\
0 \\
0 \\
\end{array}
\right).
\end{align}
Solving for the first variable 
$\mathcal{Q}_{n-1}^{(\zeta)}(z)$,
we obtain
\begin{align}
\label{quadratic ttrr for Q-n}
(u_nz+v_n)\mathcal{Q}_{n+1}^{(\zeta)}(z)=
(u_nz^2+s_nz-t_n)\mathcal{Q}_{n}^{(\zeta)}(z)-
\lambda_{n}(u_{n+1}z+v_{n+1})z
\mathcal{Q}_{n-1}^{(\zeta)}(z),
\end{align}
for $n\geq1$, with
$\mathcal{Q}_{0}^{(\zeta)}(z)=1$,
$\mathcal{Q}_{1}^{(\zeta)}(z)=
z+\sigma_1+\omega_1(\zeta)$.
The parameters involved in
\eqref{quadratic ttrr for Q-n}
are given by
\begin{align*}
u_n&:=u_n(\zeta)=\omega_{n-1}(\zeta)+\lambda_{n},
\quad
s_n:=s_n(\zeta)=\omega_n^{-1}(\zeta)u_nv_{n+1}+
v_n-\omega_{n-1}(\zeta)u_{n+1},
\\
v_n&:=v_n(\zeta)=\omega_{n-1}(\zeta)
(\omega_n(\zeta)+\sigma_{n}),
\quad
t_n:=t_n(\zeta)=\omega_n^{-1}(\zeta)
\omega_{n-1}(\zeta)\sigma_{n}v_{n+1}.
\end{align*}
However, with $\omega_n(\zeta)$ defined as in
\eqref{recursive definition of omega-n}, 
it can be seen that
$u_n\zeta+v_n=0$, $n\geq2$, 
which also yields that
$z-\zeta$ is a factor of $u_nz^2+s_nz-t_n$.
Hence cancelling the common factor $z-\zeta$
from \eqref{quadratic ttrr for Q-n}
(which is possible since
$z\in\mathbb{C}\setminus\{\zeta\}$),
we obtain the simplified form
\begin{align}
\label{eqn: ttrr for Qn with polynomial coefficients}
\mathcal{Q}_{n+1}^{(\zeta)}(z)=
\left(z+\dfrac{u_{n+1}\sigma_n\omega_{n-1}(\zeta)}
	{u_n\omega_n(\zeta)}\right)
\mathcal{Q}_{n}^{(\zeta)}(z)-
\dfrac{u_{n+1}}{u_n}\lambda_{n}z
\mathcal{Q}_{n-1}^{(\zeta)}(z),
\quad n\geq2,
\end{align}
with the initial conditions 
\begin{align}
\mathcal{Q}_1^{(\zeta)}(z)
=
z+\sigma_1+\omega_1(\zeta)
\quad\mbox{and}\quad
\mathcal{Q}_2^{(\zeta)}(z)
=
(z+\sigma_2+\omega_2(\zeta))(z+\sigma_1)-\lambda_2z,
\end{align}
found by direct computations
from \eqref{Qn as a linear combination of Bn}.
\end{proof}
Let us choose $\omega_1(\zeta)=-(\sigma_1+\zeta)$. 
With this initial choice $\omega_n(\zeta)$, $n\geq2$, 
can be uniquely generated from 
\eqref{recursive definition of omega-n}.
Further, $\mathcal{Q}_1^{(\zeta)}(z)=z-\zeta$ 
and 
$\mathcal{Q}_2^{(\zeta)}(z)$ vanishes at $z=\zeta$.
Writing $\mathcal{Q}_2^{(\zeta)}(z)=(z-\zeta)(z+y)$, 
we find that 
$y=\sigma_1+\sigma_2+\lambda_2+\omega_2(\zeta)-\zeta$
and
$-\zeta y=\sigma_1(\sigma_2+\omega_2(\zeta))$.
Thus we have
\begin{align*}
\mathcal{Q}_2^{(\zeta)}(z)=
(z-\zeta)
\left(z-\frac{\sigma_1(\sigma_2+\omega_2(\zeta))}
{\zeta}\right).
\end{align*}
We now show that 
$\mathcal{Q}_i^{(\zeta)}(z)$, $i=0,1,2$, 
also satisfy a three term recurrence relation but 
different from 
\eqref{eqn: ttrr for Qn with polynomial coefficients}.
There are various ways to do this. 
We choose $\omega_0(\zeta)=-\lambda_1$ so that
$u_1=0$ and then divide
\eqref{quadratic ttrr for Q-n} by $v_1=\lambda_1\zeta$.
Or, with $y_1,y_2,y_3$ to be determined, we write
\begin{align*}
\mathcal{Q}_2^{(\zeta)}(z)=
y_1(z+y_2)\mathcal{Q}_1^{(\zeta)}(z)-
y_3z(z-\zeta)\mathcal{Q}_0^{(\zeta)}(z).
\end{align*}
Canceling the factor $z-\zeta$, we find $y_1-y_3=1$ and 
$y_1y_2=-\sigma_1(\sigma_2+\omega_2(\zeta))\zeta^{-1}$.
We put $y_1=\pi_2\in\mathbb{C}\setminus{\{0\}}$ 
so that $y_3=\pi_2-1$ and $y_2=-\sigma_1(\sigma_2+\omega_2(\zeta))\zeta^{-1}x_2^{-1}$.
With this we have
\begin{align*}
\mathcal{Q}_2^{(\zeta)}(z)=
\pi_2\left(z-\frac{\sigma_1(\sigma_2+\omega(\zeta))}
{\pi_2\zeta}\right)\mathcal{Q}_1^{(\zeta)}(z)-
(\pi_2-1)z(z-\zeta)\mathcal{Q}_0^{(\zeta)}(z).
\end{align*}
\begin{remark}
Since $\mathcal{Q}_1^{(\zeta)}(z)=
\mathcal{Q}_2^{(\zeta)}(z)=0$ at $z=\zeta$,
from \eqref{eqn: ttrr for Qn with polynomial coefficients} 
we find that $\mathcal{Q}_n^{(\zeta)}(z)=0$ at 
the excluded point $z=\zeta$ for $n\geq1$. 
This construction depends on the unique choice of $\omega_1(\zeta)$ which makes both 
$\mathcal{Q}_1^{(\zeta)}(z)$ and
$\mathcal{Q}_2^{(\zeta)}(z)$ vanish at 
$z=\zeta$.
We would like to add that such mixed recurrence relations
were also obtained for $\mathcal{Q}_n^{(\zeta)}(z)$ in
\cite{KKB-Swami-RI-Calcolo-2018}, though for $\zeta=1$. 
\end{remark}
This motivates us to define a new polynomial sequence 
$\{\mathcal{P}_n^{(\zeta)}(z)\}_{n=0}^{\infty}$, where
\begin{align*}
\mathcal{P}_n^{(\zeta)}(z)=
\frac{\mathcal{Q}_{n+1}^{(\zeta)}(z)}{z-\zeta}=
\frac{\mathcal{B}_{n+1}(z)+
\omega_{n+1}(\zeta)\mathcal{B}_{n}(z)}{z-\zeta},
\quad n\geq0,
\end{align*}
so that $\mathcal{P}_0^{(\zeta)}(z)=1$ and
$\mathcal{P}_1^{(\zeta)}(z)=
z-\sigma_1(\sigma_2+\omega_2(\zeta))\zeta^{-1}$.
Then, \eqref{eqn: ttrr for Qn with polynomial coefficients}
gives
\begin{align}
\label{eqn: ttrr for Pn after dividing by z-1}
\mathcal{P}_{n+1}^{(\zeta)}(z)=
\left(z+\frac{u_{n+2}\sigma_{n+1}\omega_n(\zeta)}
{u_{n+1}\omega_{n+1}(\zeta)}\right)
\mathcal{P}_n^{(\zeta)}(z)-
\lambda_{n+1}\frac{u_{n+2}}{u_{n+1}}z
\mathcal{P}_{n-1}^{(\zeta)}(z),
\quad n\geq1,
\end{align}
with the initial conditions as mentioned above.
The form \eqref{eqn: ttrr for Pn after dividing by z-1} is more convenient to work with since it removes the ambiguity associated with the recurrence relation involving 
$\mathcal{Q}_i^{(\zeta)}(z)$, $i=0,1,2$, as observed above.

The next result is crucial in the sense that it helps us 
introduce invariant polynomials into our analysis.
However, as is the case, we need to impose conditions on
the parameters $\{\sigma_n\}$ and $\{\lambda_n\}$
appearing in \eqref{eqn: ttrr for Bn}. 
We do so by first choosing $\sigma_n$ such that $|\sigma_1\cdots\sigma_k|>1$ for $k\geq1$. 
Thereafter, we choose $\lambda_{n}$ such that 
\begin{align*}
	\frac{\lambda_{n+1}}{\sigma_n}=1-\left|
	\frac{\zeta}{\sigma_1\cdots\sigma_n}\right|^2,
	\quad
	n\geq1.
\end{align*}
We also restrict $\zeta\in\mathbb{C}$ such that $|\zeta|=1$.
\begin{theorem}
	\label{theorem: ranga ttrr for P-n-zeta}
Define the sequences 
$\{\hat{\tau}_n(\zeta)\}_{n=0}^{\infty}$ and
$\{\hat{\alpha}_{n-1}\}_{n=1}^{\infty}$ where
$\tau_0(\zeta)=1$, 
\begin{align}
\label{eqn: definition of tau-n and alpha-n}
\hat{\tau}_n(\zeta)=
\frac{u_{n+1}}{\omega_n(\zeta)}\sigma_1\cdots\sigma_{n}
\quad\mbox{and}\quad
\hat{\alpha}_{n-1}=-
\frac{u_{n+1}}{\hat{\tau}_n(\zeta)\omega_n(\zeta)},
\quad
n\geq1.
\end{align}
Then the polynomials 
$\mathcal{P}_{n}^{(\zeta)}(z)$, $n\geq1$,
satisfy a three term recurrence relation of the form
\begin{align}
\label{eqn: ttrr for Ranga Pnz}
\mathcal{P}_{n+1}^{(\zeta)}(z)=
(z+b_{n+1}(\zeta))\mathcal{P}_{n}^{(\zeta)}(z)-
a_{n+1}(\zeta)\mathcal{P}_{n-1}^{(\zeta)}(z), 
\quad
n\geq1,
\end{align}
with $\mathcal{P}_{0}^{(\zeta)}(z)=1$, 
$\mathcal{P}_{1}^{(\zeta)}(z)=z+b_1(\zeta)$ and where
\begin{align*}
b_n(\zeta)=\frac{\hat{\tau}_n(\zeta)}{\hat{\tau}_{n-1}(\zeta)}
\quad\mbox{and}\quad
a_{n+1}(\zeta)=[1+\hat{\tau}_n(\zeta)\hat{\alpha}_{n-1}]
[1-\overline{\zeta\hat{\tau}_n(\zeta)\hat{\alpha}_n}]\zeta,
\quad
n\geq1.
\end{align*}
\end{theorem}
\begin{proof}
With the expression for $u_{n+1}$ obtained in 
Lemma~\ref{lemma for mixed recurrence relations},
we can also write
\begin{align*}
\hat{\tau}_n(\zeta)=-
\frac{\omega_{n+1}(\zeta)+\sigma_{n+1}}{\zeta}
\sigma_1\cdots\sigma_{n}
\quad\mbox{and}\quad
\hat{\alpha}_{n-1}=-\frac{{1}}{\sigma_1\cdots\sigma_{n}},
\quad
n\geq1.
\end{align*}
This gives $\lambda_{n+1}=\sigma_n(1-|\hat{\alpha}_{n-1}|^2)$,
$n\geq1$. Further, from 
\eqref{eqn: definition of tau-n and alpha-n}
\begin{align*}
\hat{\tau}_n(\zeta)\hat{\alpha}_{n-1}=-
\frac{u_{n+1}}{\omega_n(\zeta)}
\implies
1+\hat{\tau}_n(\zeta)\hat{\alpha}_{n-1}=
1-\frac{\omega_n(\zeta)+\lambda_{n+1}}{\omega_n(\zeta)}=
-\frac{\lambda_{n+1}}{\omega_n(\zeta)},
\end{align*}
while with similar computations
\begin{align*}
1-\zeta\hat{\tau}_{n-1}(\zeta)\hat{\alpha}_{n-1}=
1+\frac{\zeta}{\sigma_n}
\left(
1+\frac{\lambda_n}{\omega_{n-1}(\zeta)}
\right)=
-\frac{\omega_n(\zeta)}{\sigma_n}
\end{align*}
Thus, we arrive at 
\begin{align*}
[1+\hat{\tau}_n(\zeta)\hat{\alpha}_{n-1}]
[1-\zeta\hat{\tau}_{n-1}(\zeta)\hat{\alpha}_{n-1}]=
\frac{\lambda_{n+1}}{\sigma_n}=
1-|\hat{\alpha}_{n-1}|^2,
\quad
n\geq1,
\end{align*}
which can be simplified to give the relation
\begin{align}
\label{eqn: relation between tau-n+1 and tau-n} 
\hat{\tau}_{n+1}(\zeta)=
\frac{\zeta\hat{\tau}_n(\zeta)-\bar{\hat{\alpha}}_n}
{1-\zeta\hat{\tau}_n(\zeta)\hat{\alpha}_n},
\quad
n\geq0,
\quad\mbox{with}\quad
\hat{\tau}_0(\zeta)=1.
\end{align}
Since $|\zeta|=1$ and $|\hat{\alpha}_n|<1$, we have $|\hat{\tau}_n|=1$ for $n\geq0$.
The following relation
\begin{align}
\label{eqn: relation in theorem-proof}
[1+\hat{\tau}_n(\zeta)\hat{\alpha}_{n-1}]
[1-\zeta\hat{\tau}_n(\zeta)\hat{\alpha}_{n}]=
\frac{\lambda_{n+1}}{\sigma_{n+1}}
\frac{\omega_{n+1}(\zeta)}{\omega_n(\zeta)},
\quad
n\geq1,
\end{align}
can also be easily derived, on either side of which 
we multiply 
\begin{align*}
\frac{\hat{\tau}_{n+1}(\zeta)}{\hat{\tau}_{n}(\zeta)}=
\frac{\omega_n(\zeta)}{\omega_{n+1}(\zeta)}
\frac{u_{n+2}}{u_{n+1}}\sigma_{n+1},
\quad
n\geq1.
\end{align*}
The right hand side of
\eqref{eqn: relation in theorem-proof}
is $\frac{\lambda_{n+1}u_{n+2}}{u_{n+1}}$.
Further, from 
\eqref{eqn: relation between tau-n+1 and tau-n},
we have
\begin{align*}
\frac{\hat{\tau}_{n+1}(\zeta)}{\hat{\tau}_n(\zeta)}
[1-\zeta\hat{\tau}_n(\zeta)\hat{\alpha}_n]=
[1-\overline{\zeta\hat{\tau}_n(\zeta)\hat{\alpha}_n}]\zeta,
\end{align*}
which used in the left hand side of 
\eqref{eqn: relation in theorem-proof}
gives
\begin{align*}
\frac{\lambda_{n+1}u_{n+2}}{u_{n+1}}=
[1+\hat{\tau}_n(\zeta)\hat{\alpha}_{n-1}]
[1-\overline{\zeta\hat{\tau}_n(\zeta)\hat{\alpha}_n}],
\quad
n\geq1.
\end{align*}
The proof is complete by comparing the above expressions with
\eqref{eqn: ttrr for Pn after dividing by z-1}.
\end{proof}
That $\mathcal{P}_n^{(\zeta)}(z)$ is self-inversive 
and satisfies 
$\hat{\tau}_n\mathcal{P}_n^{(\zeta)\ast}(z)=
\mathcal{P}_n^{(\zeta)}(z)$
can be seen by inverting the relation
\eqref{eqn: ttrr for Ranga Pnz} and using 
\eqref{eqn: relation between tau-n+1 and tau-n}. 
Further, the recurrence relation \eqref{eqn: ttrr for Ranga Pnz}
appears in the theory of para-orthogonal polynomials that are obtained from the Christoffel-Darboux (CD) kernels $K_n(w;z)$
for $|w|=1$ 
and is satisfied by the monic form of $K_n(w;z)$
\cite{Ranga-OPUC-chain-seq-JAT-2013,CMV-measures-POP-2002}.
We refer the reader to
\cite{Simon-CD-survey} for a recent survey of applications of the CD kernels in the spectral theory of 
orthogonal polynomials.
\section{A special case: the reversed Szeg\H{o} polynomials}
\label{sec: special case}
Our goal in this section is to identify a class of polynomials that satisfies
the recurrence relation \eqref{eqn: ttrr for Bn} along with the 
conditions that we have imposed on $\sigma_n$ and
$\lambda_{n+1}$, $n\geq1$.
To begin with, we are concerned only with the recurrence relation 
\eqref{eqn: ttrr for Ranga Pnz} 
and not how we arrived at it.
Then, we use the theory that associates with 
\eqref{eqn: ttrr for Ranga Pnz}, a non trivial measure on the unit circle and the corresponding orthogonal polynomials.

Restricting the value of $\zeta$ to be $1$, 
a key role is played by the scaled polynomials 
\begin{align*}
R_n(z)=\frac{\prod_{j=0}^{n-1}
[1-\hat{\tau}_j\hat{\alpha}_j]}
{\prod_{j=0}^{n-1}[1-\mathrm{Re}(\hat{\tau}_j\hat{\alpha}_j)]}
\mathcal{P}_n^{(1)}(z),
\quad
\hat{\tau}_n:=\hat{\tau}_{n}(1),
\quad
n\geq1,
\end{align*}
which, among other things, relates to a continued fraction transformation observed by Wall 
\cite[Section 78]{Wall-book}
leading to the recurrence relation
\begin{align}
\label{eqn: ttrr for R-n with c-n and d-n}
R_{n+1}(z)=[(1+ic_{n+1})z+(1-ic_{n+1})]R_{n}(z)-
4d_{n+1}zR_{n-1}(z),
\quad
n\geq1,
\end{align}
with $R_0(z)=1$ and $R_1(z)=(1+ic_1)z+(1-ic_1)$.
Here $\{c_n\}_{n=1}^{\infty}$ and $\{d_{n+1}\}_{n=1}^{\infty}$
are real sequences expressed in terms of $\hat{\tau}_n$ and 
$\hat{\alpha}_{n-1}$.
In addition, when $d_n:=(1-m_{n-1})m_n$,
$n\geq1$ ($d_n$ in such a case is called a chain sequence with 
$m_n$, $n\geq0$, the minimal parameters if $m_0=0$),
a Favard type theorem is also established 
\cite[Theorem 4.1]{Ranga-Favard-type=JAT-2014}
for 
\eqref{eqn: ttrr for R-n with c-n and d-n}
which proves the existence of a non-trivial probability measure $\hat{\mu}$ on the unit circle such that
\begin{align*}
\int_{\partial\mathbb{D}}z^{-n+k}R_n(z)(1-z)d\hat{\mu}(z)=0,
\quad
k=0,1,\cdots,n-1,
\quad
n\geq1.
\end{align*}
Further, the polynomials $\hat{S}_n(z)$ defined by
\begin{align}
\label{eqn: relation between Szego polynomials from Rn}
\hat{S}_0(z)=1,
\quad
\hat{S}_n(z)=
\frac{R_n(z)-2(1-m_n)R_{n-1}(z)}{\prod_{k=1}^{n}(1+ic_k)},
\quad
n\geq1,
\end{align}
are the monic orthogonal polynomials on the unit circle with respect to the measure $\hat{\mu}$
\cite[Theorem 5.2]{Ranga-Favard-type=JAT-2014}, 
where the associated Verblunsky coefficients 
are given by 
\begin{align*}
-\overline{\hat{S}_n(0)}=
\frac{1-2m_n-ic_n}{1+ic_n}
\prod_{k=1}^{n}\frac{1+ic_k}{1-ic_k}
\quad
n\geq1.
\end{align*}
\begin{lemma}
[\rm{\cite{Ranga-OPUC-chain-seq-JAT-2013}}]
Let us choose the real sequence $\{c_n\}_{n=1}^{\infty}$
and the positive chain sequence $\{d_n\}_{n=1}^{\infty}$ 
having the minimal parameter sequence $\{m_n\}_{n=0}^{\infty}$
where
\begin{align*}
c_k
=\frac{-\mathrm{Im}(\hat{\tau}_{k-1}\hat{\alpha}_{k-1})}
{1-\mathrm{Re}(\hat{\tau}_{k-1}\hat{\alpha}_{k-1})}
\quad\mbox{and}\quad
m_0=0,
\quad
m_k=\frac{1}{2}
\frac{|1-\hat{\tau}_{k-1}\hat{\alpha}_{k-1}|^2}
{1-\mathrm{Re}(\hat{\tau}_{k-1}\hat{\alpha}_{k-1})},
\quad
k\geq1,
\end{align*}
where $\{\hat{\tau}_k\}$ and $\{\hat{\alpha}_{k-1}\}$
are as defined in 
\eqref{eqn: definition of tau-n and alpha-n}.
Then, $\hat{\alpha}_{n-1}=-\overline{\hat{S}_n(0)}$, $n\geq1$.
\end{lemma}
\begin{proof}
We observe that with our choice for $\sigma_n$, 
$|\hat{\alpha}_{n-1}|<1$, $n\geq1$, 
and hence $\{\hat{\alpha}_{n-1}\}_{n=1}^{\infty}$ 
is eligible to
constitute a sequence of Verblunsky coefficients
\cite[Theorem 3.1.3]{Simon-book-Part-1}. 
The rest of the proof follows from simple computations using
\eqref{eqn: relation between tau-n+1 and tau-n}
for $\zeta=1$.
\end{proof}
With $\hat{\tau}_n=\prod_{j=1}^{n}\frac{1-ic_j}{1+ic_j}$, 
the following relation
\begin{align}
\label{eqn: relation between Pn and Szego polynomials}
\mathcal{P}_n^{(1)}(z)
=
\frac{1}{z-1}(z\hat{S}_n(z)-\hat{\tau}_n\hat{S}_n^{\ast}(z)),
\quad
n\geq1,
\end{align}
can also be derived from 
$\hat{S}_n(z)$ and $\hat{S}_n^{\ast}(z)$ given by
\eqref{eqn: relation between Szego polynomials from Rn},
using which it can be shown that 
$\mathcal{P}_n^{(1)}(z)$ satisfies 
the orthogonality relations
\begin{align}
\label{eqn: orthogonality relations for kernel polynomials}
\int_{\partial\mathbb{D}}z^{-n+k}\mathcal{P}_n^{(1)}(z)
(1-z)d\hat{\mu}(z)=
\begin{cases}
	\tilde{\eta_n}\neq0, &\text{$k=-1$,}\\
	0, &\text{$0\leq j\leq n-1$,}\\
	\hat{\eta}_n\neq0, &\text{$k=n$,}
\end{cases}
\end{align}
for any $n\geq1$.
Moreover, that  
$\mathcal{P}_n^{(1)\ast}(z)=
\bar{\hat{\tau}}_n\mathcal{P}_n^{(1)}(z)$,
$n\geq1$, can also be verified from 
\eqref{eqn: relation between Pn and Szego polynomials}.
\begin{remark}
\label{remark: emphasis}
The essence of the above discussion is that starting with any polynomial sequence $\mathcal{B}_n(z)$, $n\geq0$, 
satisfying \eqref{eqn: ttrr for Bn} with $\sigma_n$ and
$\lambda_n$ appropriately chosen, 
we are able to arrive at the recurrence relation
\eqref{eqn: ttrr for Ranga Pnz}
satisfied by monic CD kernels 
associated with orthogonal polynomials on the unit circle.
Thus, with $\mathcal{P}_0^{(1)}(z)=1$, 
we have on one hand
\begin{align*}
\mathcal{P}_n^{(1)}(z)=\frac{1}{z-1}
\left(
\mathcal{B}_{n+1}(z)+\omega_{n+1}\mathcal{B}_n(z)
\right),
\quad
\omega_n:=\omega_n(1),
\quad
n\geq1,
\end{align*}
while on the other 
\begin{align*}
\mathcal{P}_n^{(1)}(z)=\frac{1}{z-1}
(z\hat{S}_{n}(z)-\hat{\tau}_{n}\hat{S}_{n}^{\ast}(z)),
\quad
n\geq1.
\end{align*}
We note that the crucial link is the expression 
\eqref{eqn: definition of tau-n and alpha-n} of
the parameters 
$\hat{\tau}_n$ and $\hat{\alpha}_{n-1}$ 
in terms of $\sigma_n$, $\lambda_n+1$ and $\omega_n$, $n\geq1$.
\end{remark}
In view of 
Remark~\ref{remark: emphasis}, the special class of polynomials that we will use are the reversed Szeg\H{o} polynomials.
For any sequence $\{S_n(z)\}_{n=0}^{\infty}$ of polynomials
orthogonal on the unit circle with respect to the non-trivial positive measure $\mu$, 
consider the Szeg\H{o} recurrences
\begin{align}
\label{eqn: Szego recurrences}
\begin{split}
S_n(z)
&=zS_{n-1}(z)-\bar{\alpha}_{n-1}S_{n-1}^{\ast}(z),
\\
S_n(z)
&=(1-|\alpha_{n-1}|^2)zS_{n-1}(z)-
\bar{\alpha}_{n-1}S_n^{\ast}(z),
\quad
n\geq1.
\end{split}
\end{align}
It can be shown that the sequence of polynomials 
$\{S_n^{\ast}(z)\}_{n=0}^{\infty}$ satisfies the recurrence relation (normalized to monic form)
\begin{align*}
\frac{S_{n+1}^{\ast}(z)}{-\alpha_n}=
\left(z+\frac{\alpha_{n-1}}{\alpha_n}\right)
\frac{S_{n}^{\ast}(z)}{-\alpha_{n-1}}-
\frac{\alpha_{n-2}}{\alpha_{n-1}}(1-|\alpha_{n-1}|^2)z
\frac{S_{n-1}^{\ast}(z)}{-\alpha_{n-2}},
\quad
n\geq1,
\end{align*}
with $S_0^{\ast}(z)=1$ and 
$\frac{S_1^{\ast}(z)}{-\alpha_0}=z-\frac{1}{\alpha_0}$.
Comparing with \eqref{eqn: ttrr for Bn}, we obtain
$\mathcal{B}_0(z)=S_0^{\ast}(z)=1$,
\begin{align*}
\mathcal{B}_n(z)=\frac{S_{n}^{\ast}(z)}{-\alpha_{n-1}},
\quad
\sigma_{n+1}=\frac{\alpha_{n-1}}{\alpha_n}
\quad\mbox{and}\quad
\lambda_{n+1}=
\frac{\alpha_{n-2}}{\alpha_{n-1}}(1-|\alpha_{n-1}|^2),
\quad
n\geq1,
\end{align*}
with $|\sigma_1\cdots\sigma_{k}|=
|\frac{1}{\alpha_{k-1}}|>1$ and 
$\lambda_{k+1}=
\sigma_{k}(1-\frac{1}{|\sigma_1\cdots\sigma_k|^2})$,
$k\geq1$.
Hence, we have the following expression
\begin{align}
\label{eqn: Pn as linear combination of S-n-*}
\mathcal{P}_n^{(1)}(z)=
\frac{1}{z-1}
\left(
\frac{S_{n+1}^{\ast}(z)}{-\alpha_{n}}+\omega_{n+1}
\frac{S_{n}^{\ast}(z)}{-\alpha_{n-1}},
\right)
\quad
n\geq1.
\end{align}
which was the crux of Section~\ref{sec: Structural relations}. 
Further, from the orthogonality relations for $S_n^{\ast}(z)$
with respect to $\mu$, it follows that
\begin{align}
\label{eqn: orthogonality relation for P-n-z from S-n-*}
\int_{\partial\mathbb{D}}z^{-n+k}\mathcal{P}_n^{(1)}(z)
d\mu(z)=0,
\quad
k=0,1\cdots,n-1,
\quad
n\geq1.
\end{align}
The relation 
\eqref{eqn: orthogonality relation for P-n-z from S-n-*}
may be compared with 
\eqref{eqn: orthogonality relations for kernel polynomials},
something that we will generalize in 
Section~\ref{sec: quasi-orthogonality}.
Moreover, in Section~\ref{sec: Illustration}
we will illustrate the case when
$\hat{\alpha}_{n-1}$ 
given by \eqref{eqn: definition of tau-n and alpha-n}
is used to generate the Szeg\H{o}
polynomials via the recurrences 
\eqref{eqn: Szego recurrences}. This will explain why the class of reversed Szeg\H{o} polynomials is special.
Before that, we present a  
kind of Szeg\H{o} relation for the 
polynomials $\mathcal{P}_n^{(1)}(z)$, $n\geq0$.
\begin{lemma}
\label{lemma: relation between Pn as Szego case}
Consider the polynomial sequence $\{\mathcal{\hat{P}}_n(z)\}_{n=0}^{\infty}$
where 
\begin{align*}
\mathcal{\hat{P}}_0(z)=1, 
\quad
\hat{\tau}_{n+1}\mathcal{\hat{P}}_n(z)=
\mathcal{P}^{(1)}_{n+1}(z)-z\mathcal{P}_n^{(1)}(z),
n\geq1. 
\end{align*}
Then the following relation
\begin{align}
\label{eqn: Szego relation for P-n-z}
\mathcal{P}_n^{(1)}(z)=(1+\hat{\tau}_n\hat{\alpha}_{n-1})
(1-\hat{\tau}_n\hat{\alpha}_n)z\mathcal{P}_{n-1}^{(1)}(z)+
\hat{\tau}_n\mathcal{\hat{P}}_n(z),
\quad
n\geq1,
\end{align}
holds.
\end{lemma}
\begin{proof}
From 
\eqref{eqn: relation between tau-n+1 and tau-n} 
for $\zeta=1$,
\eqref{eqn: relation between Pn and Szego polynomials}
and the Szeg\H{o} relations 
\eqref{eqn: Szego recurrences} satisfied by
$\hat{S}_n(z)$
we obtain
\begin{align}
\label{eqn: relation between Pn1-z and Szego poly in proof}
\mathcal{P}_n^{(1)}(z)=
\frac{1}{z-1}\frac{\hat{S}_{n+1}(z)-
	\hat{\tau}_{n+1}\hat{S}_{n+1}^{\ast}(z)}
{1+\hat{\tau}_{n+1}\hat{\alpha}_{n}},
\quad
n\geq0.
\end{align}
This gives 
$z(1+\hat{\tau}_n\hat{\alpha}_{n-1})
\mathcal{P}_{n-1}^{(1)}(z)=
\mathcal{P}_n^{(1)}(z)-
\hat{\tau}_n\hat{S}_n^{\ast}(z)$, 
which upon multiplication of the factor $(1-\hat{\tau}_n\hat{\alpha}_n)$
either side leads to 
\begin{align*}
\mathcal{P}_n^{(1)}(z)=
z(1+\hat{\tau}_n\hat{\alpha}_{n-1})
(1-\hat{\tau}_n\hat{\alpha}_n)\mathcal{P}_{n-1}^{(1)}(z)+
\hat{\tau}_n
[\hat{\alpha}_n\mathcal{P}_n^{(1)}(z)+
(1-\hat{\tau}_n\hat{\alpha}_n)\hat{S}_n^{\ast}(z)]
\end{align*}
We show that the last term above is $\hat{\tau}_n\mathcal{\hat{P}}_n(z)$.
For this, we note that $\mathcal{\hat{P}}_n(z)$
has an alternative expression from the relation
\begin{align*}
\mathcal{P}_{n+1}^{(1)}(z)-z\mathcal{P}_n^{(1)}(z)=
\hat{\tau}_{n+1}
[z\hat{\alpha}_n\mathcal{P}_n^{(1)}(z)+\hat{S}_{n+1}^{\ast}(z)],
\end{align*}
where equality follows from  
\eqref{eqn: relation between Pn and Szego polynomials}.
The claim is that the relation
$
(z-1)\hat{\alpha}_n\mathcal{P}_n^{(1)}(z)=
(1-\hat{\tau}_n\hat{\alpha}_n)
\hat{S}_n^{\ast}(z)-\hat{S}_{n+1}^{\ast}(z)
$
holds. That this is true is seen from the fact that 
\begin{align*}
(z-1)\mathcal{P}_n(z)=
zS_n(z)-\tau_nS_{n}^{\ast}(z)=
S_{n+1}(z)-\tau_{n+1}(1-\tau_n\alpha_n)S_n^{\ast}(z),
\end{align*}
in which we eliminate $S_{n+1}(z)$ with its expression
obtained from 
\eqref{eqn: relation between Pn1-z and Szego poly in proof}.
\end{proof}
We conclude this section with the observation that while
$\mathcal{P}_n^{(1)}(0)=\hat{\tau}_n$, $n\geq1$, as can be seen from
\eqref{eqn: relation between Pn and Szego polynomials}
or 
\eqref{eqn: relation between Pn1-z and Szego poly in proof},
$\hat{\tau}_{n+1}\hat{\mathcal{P}}_n(0)=
\mathcal{P}_{n+1}^{(1)}(0)=\hat{\tau}_{n+1}$ which implies
$\hat{\mathcal{P}}_n(0)=1$, $n\geq1$.
\section{Quasi-orthogonality on the unit circle}
\label{sec: quasi-orthogonality}
Motivated by the representation 
\eqref{eqn: Pn as linear combination of S-n-*} that
expresses $\mathcal{P}_n^{(1)}(z)$, $n\geq1$, as a linear combination of the reversed Szeg\H{o} polynomials, we 
present the following definition of quasi-orthogonality on 
the unit circle.
\begin{definition}
\label{def: quasi-orthogonality}
Let $\mathfrak{v}$ be a positive-definite 
linear functional and $s\in\mathbb{N}$.
Let $\{\psi_n(z)\}_{n=0}^{\infty}$ be a sequence of monic polynomials. We say 
that $\psi_n(z)$ is quasi-orthogonal of order $s$ 
with respect to $\mathfrak{v}$ if
\begin{enumerate}[(i)]
\item $\mathfrak{v}[z^{-n+k}\cdot\psi_n(z)]=0$, 
for $k=s-1,s,\cdots,n-s$
and for every $n\geq2s-1$.
\item There exists $n_0\geq 2s-2$ such that 
$\mathfrak{v}[z^{-s+1}\cdot\psi_{n_0}(z)]\neq0.$ 
\end{enumerate}
\end{definition}
We note that Definition~\ref{def: quasi-orthogonality}
weakens the conditions of quasi-orthogonality presented in
\cite{Alfaro-Moral-1994} and the fact that invariant polynomials in the $\mathcal{DG}$ class are not orthogonal to constants (see 
\eqref{eqn: orthogonality relations for kernel polynomials}) plays a fundamental role.
To see that $\mathcal{P}_n^{(1)}(z)$, $n\geq1$, satisfies our definition, suppose the moment functional $\mathfrak{v}$ is
given by the integral representation 
\begin{align*}
\mathfrak{v}[z^{-n+k}\cdot\psi_n(z)]=
c\int_{\partial\mathbb{D}}z^{-n+k}\psi_n(z)
(z-1)d\mu(z),
\quad
c\in\mathbb{R}\setminus\{0\}.
\end{align*}
From \eqref{eqn: orthogonality relation for P-n-z from S-n-*},
the monic polynomials
$\mathcal{P}_n^{(1)}(z)$, $n\geq0$, satisfy
part $(i)$ of Definition~\ref{def: quasi-orthogonality}
for $s=1$.
Further, inverting the relation 
\eqref{eqn: Pn as linear combination of S-n-*}, 
we have
\begin{align*}
\bar{\hat{\tau}}_n
\int_{\partial\mathbb{D}}\mathcal{P}_n^{(1)}(z)(z-1)
d\mu(z)=
\frac{1}{-\bar{\alpha}_{n}}
\int_{\partial\mathbb{D}}S_{n+1}(z)d\mu(z)
+
\frac{\bar{\omega}_{n+1}}{-\bar{\alpha}_{n-1}}
\int_{\partial\mathbb{D}}zS_{n}(z)d\mu(z),
\end{align*}
which, upon using the Szeg\H{o} relations
\eqref{eqn: Szego recurrences} 
for the last term above gives
\begin{align*}
\int_{\partial\mathbb{D}}\mathcal{P}_n^{(1)}(z)(z-1)
d\mu(z)=
-\frac{\tau_n\bar{\omega}_{n+1}}{c}\frac{\bar{\alpha}_n}
{\bar{\alpha}_{n-1}}
||S_n(z)||^2
\neq0.
\end{align*}
Hence, the sequence 
$\{\mathcal{P}_n^{(1)}(z)\}$ given by 
\eqref{eqn: Pn as linear combination of S-n-*}
is quasi-orthogonal of order $s=1$
with respect to $\mathfrak{v}$
(where we have put $n_0=n$).

Let $\mathfrak{u}\in\mathcal{H(DG)}$, the class of moment functionals given by the integral representation
\begin{align}
\label{eqn: integral representation for u}
\mathfrak{u}[z^{-n+k}\cdot\psi_n(z)]=
\int_{\partial\mathbb{D}}z^{-n+k}\psi_n(z)
(1-z)d\delta(z)=
\begin{cases}
	\tilde{h}_n\neq0, &\text{$k=-1$,}\\
	0, &\text{$0\leq k\leq n-1$,}\\
	\hat{h}_n\neq0, &\text{$k=n$.}
\end{cases}
\end{align}
The sequence of monic polynomials 
$\{\psi_n(z)\}_{n=0}^{\infty}$
satisfying 
\eqref{eqn: integral representation for u}
is said to be associated with $\mathfrak{u}$
if $\delta$  is a non-trivial positive measure on 
$\partial{\mathbb{D}}$. 
Assuming $\psi_n(z)$, $n\geq0$,
is quasi-orthogonal with respect to $\mathfrak{v}$,
our goal in the remainder of this section is to characterize
$\mathfrak{v}$ in terms of $\mathfrak{u}\in\mathcal{H(DG)}$.

We begin with the following orthogonality properties.
If $\mathcal{P}_n^{(1)}(z)$, $n\geq1$, 
is quasi-orthogonal with respect to $\mathfrak{v}$,
then from Lemma~\ref{lemma: relation between Pn as Szego case}  
we have
$\mathfrak{v}[z^{-n+k}\cdot\mathcal{\hat{P}}_n(z)]=0$ for
$k=s-2,\cdots, n-s-1$.
Further,
$\mathfrak{u}[z^{-n+k}\cdot\mathcal{\hat{P}}_n(z)]=0$
for $k=-1,0,\cdots,n-2$.
Hence, from 
\eqref{eqn: Szego relation for P-n-z} it follows that
\begin{align}
	\label{eqn: ortho condition Pn-z-n-1}
	\mathfrak{u}[ z^{-n-1}\cdot\mathcal{P}_n^{(1)}(z)]=
	\prod_{j=1}^{n}(1+\hat{\tau}_j\hat{\alpha}_{j-1})
	(1-\hat{\tau}_j\hat{\alpha}_j),
	\quad n\geq1.
\end{align}
\begin{remark}
We invert the relation
\eqref{eqn: Szego relation for P-n-z} 
to obtain
\begin{align*}
\mathcal{P}_{n}^{(1)}(z)=
\overline{(1+\hat{\tau}_n\hat{\alpha}_{n-1})
(1-\hat{\tau}_n\hat{\alpha}_n)}
\frac{\hat{\tau}_n}{\hat{\tau}_{n-1}}
\mathcal{P}_{n-1}^{(1)}(z)+
\mathcal{\hat{P}}_{n}^{\ast}(z),
\quad
n\geq1.
\end{align*}
Since $\mathfrak{v}$ is Hermitian
\begin{align*}
\mathfrak{v}[z^{-s+1}\cdot\mathcal{\hat{P}}_{n}^{\ast}(z)]=
\mathfrak{v}[z^{-s+1}\cdot z^n\overline{\mathcal{\hat{P}}}_{n}(z^{-1})]=
\overline{
\mathfrak{v}[z^{-n+s-1}\cdot\mathcal{\hat{P}}_{n}(z)]=
}
=0,
\end{align*}
which, from 
Lemma~\ref{lemma: relation between Pn as Szego case}, 
leads to the following relation
\begin{align*}
\mathfrak{v}[z^{-s+1}\cdot\mathcal{P}_n^{(1)}(z)]=
\prod_{j=2s-1}^{n}
\overline{
(1+\hat{\tau}_j\hat{\alpha}_{j-1})
(1-\hat{\tau}_j\hat{\alpha}_j)
}
\frac{\hat{\tau}_n}{\hat{\tau}_{2s-2}}
\mathfrak{v}[z^{-s+1}\cdot\mathcal{P}_{2s-2}^{(1)}(z)],
\end{align*}
for $n\geq2s-1$. 
If in Definition~\ref{def: quasi-orthogonality} we assume 
$n_0=2s-2$, then the above relation implies
\begin{align*}
\mathfrak{v}[z^{-s+1}\cdot\mathcal{P}_{n_0}^{(1)}(z)]\neq0
\Longleftrightarrow
\mathfrak{v}[z^{-s+1}\cdot\mathcal{P}_n^{(1)}(z)]\neq0,
\quad
\forall n\geq2s-1.
\end{align*}
It may be noted that this notion is referred to as strict
quasi-orthogonality of order $s$ in 
\cite[Definition 2.3]{Alfaro-Moral-1994}.
Further, this is also reflected 
in the illustration preceding
\eqref{eqn: integral representation for u}
in the sense that quasi-orthogonality 
of order $s=1$ holds for every $n\geq0$.
\end{remark}

Consider a polynomial 
$\mathcal{U}(z)=\hat{c}_0^{(n)}+\hat{c}_1^{(n)}z+\cdots+
\hat{c}_s^{(n)}z^s$ 
with $\hat{c}_s^{(n)}\neq0$. 
Let $\mathcal{U}_{\ast}(z^{-1})$ be obtained from 
$\mathcal{U}(z)$ by replacing $z$ with $z^{-1}$
and (for notational purposes) $\hat{c}_j^{(n)}$ with $\hat{c}_{-j}^{(n)}$.
Here, the superscript $(n)$ signifies that the coefficients are associated with a polynomial of degree $n$. Consider the following representation
\begin{align}
\label{eqn: integral representation of v-u-U}
\begin{split}
\mathfrak{v}[z^{-n+k}\cdot\mathcal{P}_n^{(1)}(z)]
&=
[z^{-1}\mathcal{U}(z)+z\mathcal{U}_{\ast}(z^{-1})]
\mathfrak{u}[z^{-n+k}\cdot\mathcal{P}_n^{(1)}(z)]
\\
&=
\int_{\partial\mathbb{D}}z^{-n+k}\mathcal{P}_n^{(1)}(z)
[z^{-1}\mathcal{U}(z)+
z\mathcal{U}_{\ast}(z^{-1})](1-z)d\hat{\mu}(z).
\end{split}
\end{align}
Then, using the orthogonality properties
\eqref{eqn: integral representation for u} 
we have
\begin{align*}
\mathfrak{v}[z^{-n+k}\cdot\mathcal{P}_n^{(1)}(z)]=
\mathfrak{u}[z^{-n+k-1}
\mathcal{U}(z)\cdot\mathcal{P}_n^{(1)}(z)]+
\mathfrak{u}[z^{-n+k+1}
\mathcal{U}_{\ast}(z^{-1})\cdot\mathcal{P}_n^{(1)}(z)]=0,
\end{align*}
for
$k=s-1,s,\cdots,n-s$
and
$
\mathfrak{v}[z^{-s+1}\cdot\mathcal{P}_{2s-2}^{(1)}(z)]\neq0.
$
Hence, $\mathcal{P}_n^{(1)}(z)$, $n\geq1$, is quasi-orthogonal of order $s$ with respect to $\mathfrak{v}$.
The following puts this in
a formal setting. 
\begin{theorem}
\label{theorem: existence and uniqueness of scalars}
Let $\mathfrak{u}\in\mathcal{H(DG)}$ and 
$\{\mathcal{P}_n^{(1)}(z)\}_{0}^{\infty}$ be the sequence of polynomials associated with 
$\mathfrak{u}$.
Then, $\{\mathcal{P}_n^{(1)}(z)\}_{0}^{\infty}$ is quasi-orthogonal of order $s$ with respect to $\mathfrak{v}$ if and only if 
\begin{align}
\label{eqn: sum-v-u-theorem}
\mathfrak{v}=\sum_{j=-s+1}^{j=s-1}c_j^{(n)}z^j
\mathfrak{u},
\end{align}
where $c_j^{(n)}$ are complex scalars.
\end{theorem}
\begin{proof}
The necessity for existence of the scalars
$c_j^{(n)}$ is proved in the discussion preceding 
Theorem~\ref{theorem: existence and uniqueness of scalars}
and follows from the fact that 
$[z^{-1}\mathcal{U}(z)+z\mathcal{U}_{\ast}(z^{-1})]$ 
can in fact be
written in the required form
\eqref{eqn: sum-v-u-theorem}.
Hence, in order to prove the sufficient part, define
the linear functional
\begin{align*}
\mathfrak{w}=
\mathfrak{v}-\sum_{j=-s+1}^{s-1}c_j^{(n)}z^j\mathfrak{u},
\quad
c_j^{(n)}\in\mathbb{C}.
\end{align*}
Since $\mathcal{P}_n^{(1)}(z)$ is quasi-orthogonal with respect to $\mathfrak{v}$ and satisfies the orthogonality conditions
\eqref{eqn: integral representation for u},
it follows that
\begin{align*}
\mathfrak{w}[z^{-n+k}\cdot\mathcal{P}_n^{(1)}(z)]=0,
\quad
k=s-1,s,\cdots, n-s,
\end{align*}
and for every $n\geq2s-1$. We now prove that 
$\mathfrak{w}[z^{-n+k}\cdot\mathcal{P}_n^{(1)}(z)]=0$, separately, 
for the values 
$k=n-s+1,n-s+2,\cdots,n$ and $k=0,1,\cdots,s-2$.
The first case holds if 
\begin{align*}
\mathfrak{v}[z^{-n+k}\cdot\mathcal{P}_n^{(1)}(z)]=
\sum_{j=-s+1}^{s-1}c_j^{(n)}
\mathfrak{u}[z^{-n+k+j}\cdot\mathcal{P}_n^{(1)}(z)],
\quad
k=n-s+1,\cdots,n,
\end{align*}
which leads to the following system of equations
\begin{align*}
\lefteqn
{
\left(
\begin{array}{c}
	\mathfrak{v}[z^{-s+1}\cdot\mathcal{P}_n^{(1)}(z)] \\
	\mathfrak{v}[z^{-s+2}\cdot\mathcal{P}_n^{(1)}(z)] \\
	\vdots \\
	\mathfrak{v}[1\cdot\mathcal{P}_n^{(1)}(z)] \\
\end{array}
\right)
}
\\&&
=
\left(
\begin{array}{cccc}
\mathfrak{u}[1\cdot\mathcal{P}_n^{(1)}(z)] &  &  &  \\
\mathfrak{u}[z\cdot\mathcal{P}_n^{(1)}(z)] & \mathfrak{u}[1\cdot\mathcal{P}_n^{(1)}(z)] &  &   \\
\vdots & \vdots & \ddots &   \\
\mathfrak{u}[z^{s-1}\cdot\mathcal{P}_n^{(1)}(z)] & \mathfrak{u}[z^{s-2}\cdot\mathcal{P}_n^{(1)}(z)] & \cdots & \mathfrak{u}[1\cdot\mathcal{P}_n^{(1)}(z)] \\
\end{array}
\right)
	\left(
	\begin{array}{c}
		c_{s-1}^{(n)} \\
		c_{s-2}^{(n)} \\
		\vdots \\
		c_{0}^{(n)} \\
	\end{array}
	\right),
\end{align*}
which has a unique solution for 
$c_{0}^{(n)},\cdots,c_{s-1}^{(n)}$
since $\mathfrak{u}[\psi_n(z)\cdot 1]\neq0$.
Similarly, the second case holds if 
$c_{-s+1}^{(n)},c_{-s+2}^{(n)},
\cdots,c_{-1}^{(n)}$
satisfy the system
\begin{align*}
\lefteqn{
\left(
\begin{array}{c}
	\mathfrak{v}[z^{-n+s-2}\cdot\mathcal{P}_n^{(1)}(z)] \\
	\mathfrak{v}[z^{-n+s-3}\cdot\mathcal{P}_n^{(1)}(z)] \\
	\vdots \\
	\mathfrak{v}[z^{-n}\cdot\mathcal{P}_n^{(1)}(z)] \\
\end{array}
\right)
}
\\&&
=
\left(
\begin{array}{cccc}
		\mathfrak{u}[z^{-n-1}\cdot\mathcal{P}_n^{(1)}(z)] &  &  &  \\
		\mathfrak{u}[z^{-n-2}\cdot\mathcal{P}_n^{(1)}(z)] & \mathfrak{u}[z^{-n-1}\cdot\mathcal{P}_n^{(1)}(z)] &  &   \\
		\vdots & \vdots & \ddots &   \\
		\mathfrak{u}[z^{-n-s+1}\cdot\mathcal{P}_n^{(1)}(z)] & \mathfrak{u}[z^{-n-s+2}\cdot\mathcal{P}_n^{(1)}(z)] & \cdots & \mathfrak{u}[z^{-n-1}\cdot\mathcal{P}_n^{(1)}(z)] \\
	\end{array}
	\right)
	\left(
	\begin{array}{c}
		c_{-s+1}^{(n)} \\
		c_{-s+2}^{(n)} \\
		\vdots \\
		c_{-1}^{(n)} \\
	\end{array}
	\right)
\end{align*}
which again has a unique solution since 
$\mathfrak{u}[z^{-n-1}\cdot\mathcal{P}_n^{(1)}(z)]\neq0$.
Thus, we have shown that 
$\mathfrak{w}[z^{-n+k}\cdot\mathcal{P}_n^{(1)}(z)]=0$ for
$k=0,1,\cdots,n$ and for every $n\geq2s-1$.
Hence, arguing as in
\cite[Proposition 2.6]{Alfaro-Moral-1994}, we obtain
$\mathfrak{w}\equiv0$ which implies
$\mathfrak{v}=\sum_{j=-s+1}^{s-1}c_j^{(n)}z^j\mathfrak{u}$.
\end{proof}
The next result provides the characterization between 
$\mathfrak{v}$ and $\mathfrak{u}\in\mathcal{H(DG)}$. 
\begin{theorem}
\label{thm: relation between v and u in terms of U}
With the conditions of Theorem~\eqref{theorem: existence and uniqueness of scalars}, 
$\{\mathcal{P}_n^{(1)}(z)\}$ is quasi-orthogonal of order $s$ with respect to $\mathfrak{v}$ 
if and only if there exists an unique
polynomial $\mathcal{U}(z)$ such that
\begin{align}
\label{eqn: measure relation}
\mathfrak{v}=[z^{-1}\mathcal{U}(z)+z\mathcal{\bar{U}}(z^{-1})]
\mathfrak{u},
\end{align}
where $\deg{\mathcal{U}(z)}=s$
and $\mathcal{\bar{U}}(z^{-1})$ is obtained from 
$\mathcal{U}(z)$ by replacing $z$ with $z^{-1}$
and $c_j$ with $\bar{c}_j$, $j=0,1,\cdots,s-1$. 
\end{theorem}
\begin{proof}
Since by Theorem~\ref{theorem: existence and uniqueness of scalars} the representation 
\eqref{eqn: sum-v-u-theorem} holds,
the existence part of the proof follows 
if we show that the functional 
$\sum_{j=-s+1}^{s-1}c_j^{(n)}z^j\mathfrak{u}$
is independent of $n$ and 
$\bar{c}_{-j}=c_j$ for $j=0,1,\cdots,s-1$.
We show that the coefficients are independent of $n$
by proving that $c_j^{(n)}=c_j^{(n-1)}$ for $j=-s+1,\cdots,s-1$.

Using the orthogonality of $\hat{\mathcal{P}}_n(z)$, we
operate the functional $\mathfrak{v}$
on
\eqref{eqn: Szego relation for P-n-z}
to get
\begin{align}
\label{eqn: operate the functional to get}
\mathfrak{v}[z^{-n+s-2}\cdot\mathcal{P}_n^{(1)}(z)]=
(1+\hat{\tau}_n\hat{\alpha}_{n-1})
(1-\hat{\tau}_n\hat{\alpha}_n)
\mathfrak{v}[z^{-n+s-1}\cdot\mathcal{P}_n^{(1)}(z)].
\end{align}
On the other hand, from the system of equations
involving $c_{-j}^{(n)}$, $j=s-1,\cdots,1$,
obtained in the proof of Theorem~\ref{theorem: existence and uniqueness of scalars}, 
we have
\begin{align*}
\mathfrak{v}[z^{-n+s-2}\cdot\mathcal{P}_n^{(1)}(z)]=
c_{-s+1}^{(n)}\mathfrak{u}
[z^{-n-1}\cdot\mathcal{P}_n^{(1)}(z)].
\end{align*}
which used with \eqref{eqn: operate the functional to get}
sets up an iteration leading to
\begin{align*}
c_{-s+1}^{(n)}\mathfrak{u}
[z^{-n-1}\cdot\mathcal{P}_n^{(1)}(z)]=
(1+\hat{\tau}_{n}\hat{\alpha}_{n-1})
(1-\hat{\tau}_n\hat{\alpha}_n)
\prod_{j=1}^{n-1}
(1+\hat{\tau}_{j}\hat{\alpha}_{j-1})
(1-\hat{\tau}_j\hat{\alpha}_j)
c_{-s+1}^{(n-1)}.
\end{align*}
That $c_{-s+1}^{(n)}=c_{-s+1}^{(n-1)}$
immediately follows from 
\eqref{eqn: ortho condition Pn-z-n-1}.
We now assume $c_{j}^{(n)}=c_{j}^{(n-1)}$
for $j=-s+1,\cdots,t-1$, where 
$t=-s+1,-s+2,\cdots,-1$.
We have
\begin{align*}
\mathfrak{v}[z^{-n-t-1}\cdot\mathcal{P}_n^{(1)}(z)]=
\sum_{j=-s+1}^{t}c_j^{(n)}\mathfrak{u}
[z^{-n-1+j-t}\cdot\mathcal{P}_n^{(1)}(z)],
\end{align*}
while Lemma~\ref{lemma: relation between Pn as Szego case}
gives $\tau_{n+1}\mathfrak{v}[z^{-n-t-1}\cdot
\mathcal{\hat{P}}_n^{(1)}(z)]$
\begin{align*}
=
\sum_{j=-s+1}^{t-1}c_j^{(n+1)}
\mathfrak{u}[z^{-n-1+j-t}\cdot\mathcal{P}_{n+1}^{(1)}(z)]-
\sum_{j=-s+1}^{t-1}c_j^{(n)}
\mathfrak{u}[z^{-n+j-t}\cdot\mathcal{P}_n^{(1)}(z)].
\end{align*}
Using the above two relations in 
\eqref{eqn: Szego relation for P-n-z} and separating 
the case $j=t$, it follows from the induction hypothesis that
\begin{align*}
	\lefteqn{
c_t^{(n)}\mathfrak{u}[z^{-n-1}\cdot\mathcal{P}_n^{(1)}(z)]-
(1+\hat{\tau}_n\hat{\alpha}_{n-1})
(1-\hat{\tau}_n\hat{\alpha}_n)
\mathfrak{u}[z^{-n}\cdot\mathcal{P}_{n-1}^{(1)}(z)]
c_t^{(n-1)}
}
\\&&=
\sum_{j=-s+1}^{t-1}
\left[
c_j^{(n-1)}
(1+\hat{\tau}_n\hat{\alpha}_{n-1})
(1-\hat{\tau}_n\hat{\alpha}_n)
\mathfrak{u}[z^{-n+j-t}\cdot\mathcal{P}_{n-1}^{(1)}(z)]
-
c_j^{(n)}\mathfrak{u}[z^{-n-1+j-t}\cdot
\mathcal{P}_{n}^{(1)}(z)]
\right]
\\&&+
\frac{\hat{\tau}_n}{\hat{\tau}_{n+1}}
\sum_{j=-s+1}^{t-1}c_j^{(n+1)}
\mathfrak{u}[z^{-n-1+j-t}\cdot\mathcal{P}_{n+1}^{(1)}(z)]
-
\sum_{j=-s+1}^{t-1}
\frac{\hat{\tau}_n}{\hat{\tau}_{n+1}}c_j^{(n)}
\mathfrak{u}[z^{-n+j-t}\cdot\mathcal{P}_{n}^{(1)}(z)].
\end{align*}
The induction is complete if the right hand side
above vanishes. This is shown using 
the three term recurrence relation
\eqref{eqn: ttrr for Ranga Pnz} 
for $\zeta=1$ and observing that
\begin{align*}
b_{n+1}=\frac{\hat{\tau}_{n+1}}{\hat{\tau}_n}
\quad\mbox{and}\quad
\frac{\hat{\tau}_n}{\hat{\tau}_{n+1}}a_{n+1}(1)=
(1+\hat{\tau}_n\hat{\alpha}_{n-1})
(1-\hat{\tau}_n\hat{\alpha}_n).
\end{align*}
Similar computations prove 
$c_j^{(n)}=c_j^{(n-1)}$ for $j=0,\cdots,s-1$.
Hence, we put
\begin{align*}
c_j^{(n)}=c_{j}^{(n-1)}=c_j,
\quad
j=-s+1,\cdots,s-1, 
\end{align*}
which we use to prove the next claim in the theorem.

The key point in this part of the proof is that  
$\mathfrak{u,v}$ are Hermitian and
$\mathcal{P}_{n}^{(1)}(z)$ is a self-inversive 
polynomial. Hence
\begin{align*}
\mathfrak{v}[z^{-n+s-2}\cdot\mathcal{P}_{n}^{(1)}(z)]=
c_{-s+1}^{(n)}\mathfrak{u}
[z^{-n-1}\cdot\mathcal{P}_{n}^{(1)}(z)]
\Longrightarrow
\mathfrak{v}[z^{-s+2}\cdot\mathcal{P}_{n}^{(1)}(z)]=
\overline{c_{-s+1}^{(n)}}\mathfrak{u}
[z\cdot\mathcal{P}_{n}^{(1)}(z)].
\end{align*}
With $\mathfrak{v}[z^{-s+1}\cdot\mathcal{P}_{n}^{(1)}(z)]
=
c_{s-1}^{(n)}\mathfrak{u}[1\cdot\mathcal{P}_{n}^{(1)}(z)]$,
the relation
\eqref{eqn: Szego relation for P-n-z} yields
\begin{align*}
\lefteqn{
c_{s-1}^{(n)}
\left[
\frac{\hat{\tau}_n}{\hat{\tau}_{n+1}}
\mathfrak{u}[1\cdot\mathcal{P}_{n+1}^{(1)}(z)]-
\mathfrak{u}[1\cdot\mathcal{P}_{n}^{(1)}(z)]
\right]
}
\\&&=
\overline{c_{-s+1}^{(n)}}
\left[
\frac{\hat{\tau}_n}{\hat{\tau}_{n+1}}
\mathfrak{u}[z\cdot\mathcal{P}_{n}^{(1)}(z)]-
(1+\hat{\tau}_n\hat{\alpha}_{n-1})
(1-\hat{\tau}_n\hat{\alpha}_n)
\mathfrak{u}[z\cdot\mathcal{P}_{n-1}^{(1)}(z)]
\right].
\end{align*}
We let $\mathfrak{u}$ operate on
\eqref{eqn: ttrr for Ranga Pnz} to conclude
that $c_{s-1}=\bar{c}_{-s+1}$.
We now assume 
$c_{j}=\bar{c}_{-j}$ for
$j=-t+1,\cdots,s-1$ where 
$t=-s+1,-s+2,\cdots,-1$.
The induction is completed by using the relations
\begin{align*}
\mathfrak{v}[z^{t}\cdot\mathcal{P}_{n}^{(1)}(z)]
&=
\sum_{j=-t}^{s-1}c_j\mathfrak{u}
[z^{j+t}\cdot\mathcal{P}_{n}^{(1)}(z)],
\\
\mathfrak{v}[z^{t+1}\cdot\mathcal{P}_{n-1}^{(1)}(z)]
&=
\sum_{j=-t}^{s-1}\bar{c}_{-j}\mathfrak{u}
[z^{j+t+1}\cdot\mathcal{P}_{n-1}^{(1)}(z)],
\end{align*}
in \eqref{eqn: Szego relation for P-n-z},
separating the case $j=-t$ and letting 
$\mathfrak{u}$ operate on 
\eqref{eqn: ttrr for Ranga Pnz}.
Finally, from 
\eqref{eqn: sum-v-u-theorem},
since $\mathfrak{v}[1]\in\mathbb{R}$, it follows that 
$c_0\in\mathbb{R}$. Thus writing
\begin{align*}
\mathcal{U}(z)=\frac{1}{2}(c_{s-1}z^s+\cdots+c_2z^3+
c_1z^2+c_0z+\bar{c}_{1}),
\end{align*}
the expression
\eqref{eqn: measure relation} follows.

To prove uniqueness of $\mathcal{U}(z)$, 
let $\mathcal{U}_1(z)$ and $\mathcal{U}_2(z)$
be any two polynomial solutions of 
\eqref{eqn: measure relation}, where 
$\deg{\mathcal{U}_1(z)}=s_1$ and 
$\deg{\mathcal{U}_2(z)}=s_2$.
Let $\mathcal{U}_3(z)=\mathcal{U}_2(z)-\mathcal{U}_1(z)$, 
with
$\deg{\mathcal{U}_3(z)}=r=\max\{s_1,s_2\}$.
Then from \eqref{eqn: measure relation}, we have
\begin{align*}
[z^{-1}\mathcal{U}_3(z)+z\bar{\mathcal{U}}_3(z^{-1})]
\mathfrak{u}=
\sum_{j=-r+1}^{r-1}u_jz^j\mathfrak{u}=0.
\end{align*}
We let the above functional act on 
$[z^{-n+k}\cdot\mathcal{P}_{n}^{(1)}(z)]$
for $n\geq2r-1$ and $k\geq-1$ to have
\begin{align*}
\sum_{j=-r+1}^{r-1}u_j
\mathfrak{u}[z^{-n+k+j}\cdot\mathcal{P}_{n}^{(1)}(z)]=0,
\end{align*}
which leads to a homogeneous system of equations for the choice
of $k=-1,0,\cdots,r-2$. 
The solution is trivial with
$u_{-r+1}=u_{-r+2}=\cdots=u_0=0$, so that 
$\mathcal{U}_1(z)=\mathcal{U}_2(z)$, 
thus establishing uniqueness.
\end{proof}
Observe that for a quasi-orthogonal polynomial of order $s=1$, $\mathcal{U}(z)=c_0/2$, so that 
$\mathfrak{v}=c_0\mathfrak{u}$. This fact
is reflected in the illustration preceding
\eqref{eqn: integral representation for u}. 
\section{Illustration}
\label{sec: Illustration}
Explicit representations of orthogonal polynomials on the unit circle are available in the literature
\cite{Simon-book-Part-1,Szego-book}. 
In this section, through a judicious use of the contiguous relations satisfied by hypergeometric functions, we will illustrate the theory presented so far.
The hypergeometric functions are denoted as
\begin{align*}
F(a,b;c;z)=
\sum_{k=0}^{\infty}\frac{(a)_k(b)_k}{(c)_k}
\frac{z^k}{k!},
\quad
|z|<1,
\end{align*}
where $a,b\in\mathbb{C}$ and
$c\in\mathbb{C}\setminus\{0\}$.
For the theory of hypergeometric functions, 
the contiguous relations satisfied by them and
the Pochhammer symbol $(a)_k$, we refer to
\cite{Andrews-Askey-Roy-1999}.

A family of Szeg\H{o} polynomials in terms of hypergeometric functions is given by
\cite{Ranga-OPUC-hypergeometric-AMS-2010}
\begin{align}
\label{eqn: class of Szego polynomials-illustration}
\begin{split}
S_n(z)
&=
\frac{(b+\bar{b}+1)_n}{(b+1)_n}F(-n,b+1;b+\bar{b+1};1-z),
\quad
n\geq1,
\\
S_n^{\ast}(z)
&=
\frac{(b+\bar{b}+1)_n}{(\bar{b}+1)_n}F(-n,b;b+\bar{b+1};1-z),
\quad
n\geq1,
\end{split}
\end{align}
which are orthogonal on the unit circle with respect to the measure 
\cite[Theorem 4.1]{Ranga-OPUC-hypergeometric-AMS-2010}
$d\mu(b;e^{i\theta})=\phi(b;\theta)d\theta$,
where for $0\leq\theta\leq2\pi$
\begin{align*}
\phi(b;\theta)=\tau^{b}e^{(\pi-\theta)\mathrm{Im}{b}}[\sin{(\theta/2)}]^{2\mathrm{Re}{b}},
\quad
\tau^{b}=\frac{2^{b+\bar{b}}|\Gamma(b+1)|^2}
{2\pi\Gamma(b+\bar{b}+1)}.
\end{align*}
We first express the Szeg\H{o} recurrences
\eqref{eqn: Szego recurrences} 
satisfied by $S_n(z)$ and $S_n^{\ast}(z)$ using 
contiguous relations.
Consider
\begin{align*}
(c-a-b)F(a,b;c;1-z)+azF(a+1,b;c;1-z)=
(c-b)F(a,b-1;c;1-z),
\end{align*}
in which we substitute $a=-n$, $b=b+1$, $c=b+\bar{b}+1$
and multiply $\frac{(b+\bar{b}+1)_n}{(b+1)_n}$ 
to get
\begin{align}
\label{eqn: contiguous relation-Szego-relation-1}
S_n(z)=\frac{n(b+\bar{b}+n)}{|b+n|^2}zS_{n-1}(z)+
\frac{(\bar{b})_n}{(b+1)_n}S_n^{\ast}(z), 
\quad
n\geq1.
\end{align}
Next, in the contiguous relation
%
\begin{align*}
F(a,b+1;c;1-z)=\frac{z(a-b)}{a-c+1}F(a+1,b+1;c;1-z)-
\frac{(c-b-1)}{(a-c+1)}F(a+1,b;c;1-z),
\end{align*}
we substitute $a=-n$, $c=b+\bar{b}+1$ and multiply
$\frac{(b+\bar{b}+1)_n}{(b+1)_n}$ both sides to get
\begin{align}
	\label{eqn: contiguous relation-Szego-relation-2}
S_n(z)=zS_{n-1}+
\frac{(\bar{b})_n}{(b+1)_n}S_{n-1}^{\ast}(z).
\end{align}
We observe that the relations
\eqref{eqn: contiguous relation-Szego-relation-1}
and \eqref{eqn: contiguous relation-Szego-relation-2}
constitute the Szeg\H{o} recurrences
\eqref{eqn: Szego recurrences} if we identify
$\bar{\alpha}_{n-1}=-\frac{(\bar{b})_n}{(b+1)_n}$ so that
$1-|\alpha_{n-1}|^2=\frac{n(b+\bar{b}+1)}{|b+n|^2}$.
The first equality follows from
\cite[Corollary 2.2.3]{Andrews-Askey-Roy-1999}
\begin{align*}
-\bar{\alpha}_{n-1}=S_n(0)=
\frac{(b+\bar{b}+1)_n}{(b+1)_n}F(-n,b+1;b+\bar{b}+1;1),
\end{align*}
while the second equality follows from simple computations.
We now illustrate the characterization 
\eqref{eqn: Pn as linear combination of S-n-*} 
as discussed in Section~\ref{sec: special case}.
Using the Szeg\H{o}
relations \eqref{eqn: contiguous relation-Szego-relation-1}
and \eqref{eqn: contiguous relation-Szego-relation-2}
along with the expressions for $\alpha_{n-1}$, 
we begin with the recurrence relation
\begin{align*}
	\frac{S_{n+1}^{\ast}(z)}{-\alpha_n}=
	\left(z+\frac{\bar{b}+n+1}{b+n}\right)
	\frac{S_{n}^{\ast}(z)}{-\alpha_{n-1}}-
	\frac{n(b+\bar{b}+n)}{(b+n-1)(b+n)}z
	\frac{S_{n-1}^{\ast}(z)}{-\alpha_{n-2}},
	\quad
	n\geq1,
\end{align*}
with $S_0^{\ast}(z)=1$ and 
$\frac{S_{1}^{\ast}(z)}{-\alpha_{0}}=z+\frac{\bar{b}+1}{b}$,
which is satisfied by the monic polynomials
\begin{align*}
	\frac{S_{n}^{\ast}(z)}{-\alpha_{n-1}}=
	\frac{(b+\bar{b}+1)_n}{(b)_n}F(-n,b;b+\bar{b}+1;1-z),
	\quad
	n\geq1.
\end{align*}
It can be easily verified from 
\eqref{recursive definition of omega-n} that
\begin{align*}
\sigma_{n}=\frac{\bar{b}+n}{b+n-1}
\hbox{ and } 
\lambda_{n+1}=\frac{n(b+\bar{b}+n)}{(b+n-1)(b+n)}
\Longrightarrow
\omega_n=-\frac{b+\bar{b}+n}{b+n-1},
\quad
n\geq1.
\end{align*}
Thus, we have
\begin{align}
\label{eqn: P-n-z as linear combination of S-n-*-illustration}
\mathcal{P}_n^{(1)}(z)=
\frac{1}{z-1}
\left(
\frac{S_{n+1}^{\ast}(z)}{-\alpha_n}-
\frac{b+\bar{b}+n+1}{b+n}
\frac{S_{n}^{\ast}(z)}{-\alpha_{n-1}}
\right),
\quad
n\geq1.
\end{align}
The sequence of polynomials $\mathcal{P}_n^{(1)}(z)$, $n\geq0$,
is thus quasi-orthogonal of order $s=1$ on the unit circle
with respect to $\phi(b;\theta)d\theta$
since using the orthogonality properties of $S_n^{\ast}(z)$  
and \cite[Theorem 3.1]{Ranga-OPUC-hypergeometric-AMS-2010}
we have
\begin{align*}
	\int_{\partial\mathbb{D}}z^{-n+k}
	\mathcal{P}_n^{(1)}(z)(1-z)\phi(b;\theta)d\theta
	=
	\begin{cases} 
		0 &\mbox{if } k=0,\cdots,n-1, \\
		(h_n^{(b)})^{-2} & \mbox{if } k=n, 
	\end{cases}
\end{align*}
where the orthogonality constant is given by
\begin{align*}
	h_n^{(b)}=||S_n(z)||^{-1}=
	\sqrt{\frac{|(b+1)_n|^2}{(b+\bar{b}+1)_n}n!},
	\quad
	n\geq1.
\end{align*}
Further, an expression for $\mathcal{P}_n^{(1)}(z)$ given by
\eqref{eqn: P-n-z as linear combination of S-n-*-illustration}
can be found from the contiguous relation
\begin{align*}
\frac{b}{c}(z-1)F(a,b+1;c+1;1-z)=F(a-1,b;c;1-z)-F(a,b;c;1-z),
\end{align*}
in which we substitute $a=-n$, $c=b+\bar{b}+1$
and multiply 
$\frac{(b+\bar{b}+1)_{n+1}}{(b)_{n+1}}$ 
both sides to get
\begin{align}
\label{eqn: expression for Pn-z in illustration}
\mathcal{P}_n^{(1)}(z)=
\frac{(b+\bar{b}+2)_n}{(b+1)_n}
F(-n,b+1,b+\bar{b}+2;1-z),
\quad
n\geq1.
\end{align}
We now find the Szeg\H{o} polynomials 
$\hat{S}_n(z)$, $n\geq1$, which appear in 
\eqref{eqn: relation between Pn and Szego polynomials}
or
\eqref{eqn: relation between Pn1-z and Szego poly in proof}
by first finding the scaled polynomials $R_n(z)$.
From the definitions \eqref{eqn: definition of tau-n and alpha-n} we have
\begin{align*}
\hat{\tau}_n
=
-(\omega_{n+1}+\sigma_{n+1})\sigma_1\cdots\sigma_{n}=
\frac{(\bar{b}+1)_n}{(b+1)_n}
\hbox{ and }
\hat{\alpha}_{n-1}=-\frac{1}{\sigma_1\cdots\sigma_{n}}
=-\frac{(b)_n}{(\bar{b}+1)_n}
\end{align*}
which are the same as obtained in 
\eqref{eqn: class of Szego polynomials-illustration}
earlier.
With $\lambda=\mathrm{Re}(b)$, 
the scaled polynomials are given by
\cite{Ranga-Favard-type=JAT-2014,Ranga-OPUC-chain-seq-JAT-2013}
\begin{align*}
R_n(z)=\frac{(2\lambda+2)_n}{(\lambda+1)_n}
F(-n,b+1;b+\bar{b}+2;1-z),
\quad
n\geq1,
\end{align*}
and the corresponding Szeg\H{o} polynomials $\hat{S}_n(z)$
obtained in 
\eqref{eqn: relation between Szego polynomials from Rn}
are the same as given by
\eqref{eqn: class of Szego polynomials-illustration}.
Further, with
$1+\hat{\tau}_n\hat{\alpha}_{n-1}=\frac{n}{b+n}$,
$n\geq1$,
%
%
%
we find an expression for 
\eqref{eqn: relation between Pn1-z and Szego poly in proof},
for instance, using the contiguous relation
\begin{align*}
F(a,b;c;1-z)-F(a,b+1;c;1-z)=
\frac{a}{c}(z-1)F(a+1,b+1;c+1;1-z).
\end{align*}
We substitute $a=-n$, $c=b+\bar{b}+1$ and multiply
$\frac{(b+\bar{b}+1)_n}{(b+1)_n}$ both sides to get
\begin{align*}
S_n(z)-\frac{(\bar{b}+1)_n}{(b+1)_n}S_n^{\ast}(z)=
\frac{n}{b+n}(z-1)
\frac{(b+\bar{b}+2)_{n-1}}{(b+1)_{n-1}}
F(-n+1,b+1;b+\bar{b}+2;1-z),
\end{align*}
which yields an expression for $\mathcal{P}_n^{(1)}(z)$ as
obtained in
\eqref{eqn: expression for Pn-z in illustration}.
The expression 
\eqref{eqn: relation between Pn and Szego polynomials}
can also be found on similar lines from the contiguous relation
\begin{align*}
c(1-z)F(a,b;c;z)-cF(a-1,b;c;z)+(c-b)zF(a,b;c+1;z)=0.
\end{align*}
We conclude with the final remark that to begin with the class
of reversed Szeg\H{o} polynomials is a special case because
it allowed us to express members of the 
$\mathcal{DG}$ class of invariant polynomials in two different ways using the same sequence of Szeg\H{o} polynomials.
As might have been observed, the key feature we have used is that these invariant polynomials satisfy a three term recurrence relation, which is 
also a key feature of quasi-orthogonal polynomials on the real line.


\begin{thebibliography}{99}
%
\bibitem{Alfaro-Moral-1994}
M. Alfaro\ and\ L. Moral, 
Quasi-orthogonality on the unit circle 
and semi-classical forms, 
Portugal. Math. {\bf 51} (1994), no.~1, 47--62.
%
\bibitem{Andrews-Askey-Roy-1999}
G. E. Andrews, R. Askey\ and\ R. Roy, {\it Special functions},
Encyclopedia of Mathematics and its Applications, 71,
Cambridge Univ. Press, Cambridge, 1999.
%
\bibitem{KKB-Swami-RI-Calcolo-2018}
K. K. Behera\ and\ A. Swaminathan,
Biorthogonality and para-orthogonality of $R_I$ polynomials,
Calcolo
{\bf 55} (2018), no.~4, Art. 41, 22 pp.
%
%
\bibitem{Branquinho-Marcellan-QOP-1996}
A. Branquinho\ and\ F. Marcell\'{a}n, Generating new classes of orthogonal polynomials, Internat. J. Math. Math. Sci. {\bf 19} (1996), no.~4, 643--656.
%
\bibitem{Breuer-POP-singular-JAT-2020}
J. Breuer\ and\ E. Seelig, On the spacing of zeros of paraorthogonal polynomials for singular measures, J. Approx. Theory {\bf 259} (2020), 105482, 20 pp.
%
\bibitem{Brezinski-Driver-QOP-2004}
C. Brezinski, K. A. Driver\ and\ M. Redivo-Zaglia, Quasi-orthogonality with applications to some families of classical orthogonal polynomials, Appl. Numer. Math. {\bf 48} (2004), no.~2, 157--168.
%
\bibitem{Bultheel-quasi-POP_JCAM-2022}
A. Bultheel, R. Cruz-Barroso\ and\ C. D\'{\i}az Mendoza, 
Zeros of quasi-paraorthogonal polynomials and positive quadrature, J. Comput. Appl. Math. {\bf 407} (2022), Paper No. 114039.
%
\bibitem{CMV-measures-POP-2002}
M. J. Cantero, L. Moral\ and\ L. Vel\'{a}zquez, 
Measures and para-orthogonal polynomials 
on the unit circle, 
East J. Approx. {\bf 8} (2002), no.~4, 447--464.
%
\bibitem{Ranga-Favard-type=JAT-2014}
K. Castillo, M.S. Costa, A. Sri Ranga\ and\ D.O. Veronese,
A Favard type theorem for orthogonal polynomials on the unit
circle from a three term recurrence formula,
J. Approx. Theory
{\bf 184} (2014), 146--162.
%
\bibitem{Chihara-QOP-PAMS-1957}
T. S. Chihara, 
On quasi-orthogonal polynomials, 
Proc. Amer. Math. Soc. {\bf 8} (1957), 765--767.
%
\bibitem{Ranga-OPUC-chain-seq-JAT-2013}
M. S. Costa, H. M. Felix\ and\ A. Sri Ranga,
Orthogonal polynomials on the unit circle and
chain sequences,
J. Approx. Theory
{\bf 173} (2013), 14--32.
%
\bibitem{DG-split-Levinson-1986}
P. Delsarte\ and\ Y. V. Genin, The split Levinson algorithm, IEEE Trans. Acoust. Speech Signal Process. {\bf 34} (1986), no.~3, 470--478.
%
\bibitem{DG-tridiagonal-Szego-1988}
P. Delsarte\ and\ Y. Genin, The tridiagonal approach to Szeg\H{o}'s orthogonal polynomials, Toeplitz linear systems, and related interpolation problems, SIAM J. Math. Anal. {\bf 19} (1988), no.~3, 718--735.
%
\bibitem{DG-OPUC-signal-1989}
P. Delsarte\ and\ Y. Genin, On the role of orthogonal polynomials on the unit circle in digital signal processing applications, in {\it Orthogonal polynomials (Columbus, OH, 1989)}, 115--133, NATO Adv. Sci. Inst. Ser. C: Math. Phys. Sci., 294, Kluwer Acad. Publ., Dordrecht.
%
\bibitem{DG-tridiagonal-algebraic-I-1991}
P. Delsarte\ and\ Y. Genin, Tridiagonal approach to the algebraic environment of Toeplitz matrices. I. Basic results, SIAM J. Matrix Anal. Appl. {\bf 12} (1991), no.~2, 220--238. 
%
\bibitem{DG-trdiagonal-algebraic-II-1991}
P. Delsarte\ and\ Y. Genin, Tridiagonal approach to the algebraic environment of Toeplitz matrices. II. Zero and eigenvalue problems, SIAM J. Matrix Anal. Appl. {\bf 12} (1991), no.~3, 432--448.
%
\bibitem{Dewilde-Dym-lossless-prediction-1981}
P. Dewilde\ and\ H. Dym, Lossless chain scattering matrices and optimum linear prediction: the vector case, Internat. J. Circuit Theory Appl. {\bf 9} (1981), no.~2, 135--175. 
%
\bibitem{Dickinson-QOP-PAMS-1961}
D. Dickinson, 
On quasi-orthogonal polynomials, 
Proc. Amer. Math. Soc. 12 (1961) 185–194.
%
\bibitem{Draux-quasi-orthogonality-JAT-1990}
A. Draux,
On quasi-orthogonal polynomials,
J. Approx. Theory
{\bf 62} (1990), no.~1, 1--14.%
%
\bibitem{Golinskii-POP-2002}
L. Golinskii, 
Quadrature formula and zeros of para-orthogonal 
polynomials on the unit circle, 
Acta Math. Hungar. 
{\bf 96} (2002), no.~3, 169--186.
%
\bibitem{Grenander-Szego-book-Toeplitz-form}
U. Grenander\ and\ G. Szeg\H{o}, 
{\it Toeplitz forms and their applications}, 
second edition, Chelsea Publishing Co., New York, 1984.
%
\bibitem{Ismail-Wang-QOP-ODE-2019}
M. E. H. Ismail\ and\ X.-S. Wang, 
On quasi-orthogonal polynomials: 
their differential equations, discriminants and electrostatics, J. Math. Anal. Appl. {\bf 474} (2019), no.~2, 1178--1197. 
%
\bibitem{Jones-Njastad-Szego-digital-1991}
W. B. Jones\ and\ O. Nj\aa stad, 
Applications of Szeg\H{o} polynomials to digital signal processing, 
Rocky Mountain J. Math. {\bf 21} (1991), no.~1, 387--436.
%
\bibitem{Jones-Njastad-Thron-moment-BLMS-1989}
W. B. Jones, O. Nj\aa stad\ and\ W. J. Thron,
Moment theory, orthogonal polynomials, quadrature, and
continued fractions associated with the unit circle,
Bull. London Math. Soc.
{\bf 21} (1989), no.~2, 113--152.
%
\bibitem{Jones-Njastad-Waadeland-Szego-frequency-1994}
W. B. Jones, O. Nj\aa stad\ and\ H. Waadeland, 
Application of Szeg\H{o} polynomials to frequency analysis, SIAM J. Math. Anal. 
{\bf 25} (1994), no.~2, 491--512.
%
\bibitem{Joulak-QOP-ANM-2005}
H. Joulak, 
A contribution to quasi-orthogonal 
polynomials and associated polynomials, 
Appl. Numer. Math. 54 (2005)
65–78.
%
\bibitem{Marcellan-QOP-Berstein-1996}
F. Marcell\'{a}n, F. Peherstorfer\ and\ R. Steinbauer, Orthogonality properties of linear combinations of orthogonal polynomials, Adv. Comput. Math. {\bf 5} (1996), no.~4, 281--295.
%
\bibitem{Shohat-QOP-1937}
J. Shohat, On mechanical quadratures, 
in particular, with positive coefficients, 
Trans. Amer. Math. Soc. 
{\bf 42} (1937), no.~3, 461--496.
%
\bibitem{Simon-book-Part-1}
B. Simon,
{\it Orthogonal polynomials on the unit circle. Part 1},
American Mathematical Society Colloquium Publications, 54, Part 1, American Mathematical Society, Providence, RI, 2005.
%
\bibitem{Simon-CD-survey}
B. Simon, The Christoffel-Darboux kernel, in {\it Perspectives in partial differential equations, harmonic analysis and applications}, 295--335, Proc. Sympos. Pure Math., 79, Amer. Math. Soc., Providence, RI.
%
\bibitem{Ranga-OPUC-hypergeometric-AMS-2010}
A. Sri Ranga,
Szeg\H{o} polynomials from hypergeometric functions,
Proc. Amer. Math. Soc.
{\bf 138} (2010), no.~12, 4259--4270.
%
%
\bibitem{Szego-book}
G. Szeg\H{o}, 
{\it Orthogonal polynomials}, 
fourth edition, American Mathematical Society Colloquium Publications, Vol. XXIII, 
American Mathematical Society, Providence, RI, 1975.
%
\bibitem{Wall-book}
H. S. Wall, {\it Analytic Theory of Continued Fractions}, D. Van Nostrand Co., Inc., New York, NY, 1948.
%
\bibitem{Xu-QOP-JMAA-1994}
Y. Xu, Quasi-orthogonal polynomials, quadrature, and interpolation, J. Math. Anal. Appl. {\bf 182} (1994), no.~3, 779--799.
%
\end{thebibliography}
\end{document}